\documentstyle[amssymb,12pt]{amsart}

\newtheorem{lem}{Lemma}
\newtheorem{prop}{Proposition}

\theoremstyle{definition}

\theoremstyle{definition}
\newtheorem{thm}{Theorem}

\theoremstyle{remark}
\newtheorem{rem}{Remark}

\numberwithin{equation}{section}

\begin{document}

\newcommand{\thmref}[1]{Theorem~\ref{#1}}
\newcommand{\secref}[1]{Sect.~\ref{#1}}
\newcommand{\lemref}[1]{Lemma~\ref{#1}}
\newcommand{\propref}[1]{Proposition~\ref{#1}}
\newcommand{\corref}[1]{Corollary~\ref{#1}}
\newcommand{\remref}[1]{Remark~\ref{#1}}
\newcommand{\nc}{\newcommand}
\nc{\on}{\operatorname}
\nc{\ch}{\mbox{ch}}
\nc{\Z}{{\Bbb Z}}
\nc{\C}{{\Bbb C}}
\nc{\pone}{{\Bbb C}{\Bbb P}^1}
\nc{\pa}{\partial}
\nc{\F}{{\cal F}}
\nc{\arr}{\rightarrow}
\nc{\larr}{\longrightarrow}
\nc{\al}{\alpha}
\nc{\ri}{\rangle}
\nc{\lef}{\langle}
\nc{\W}{{\cal W}}
\nc{\la}{\lambda}
\nc{\ep}{\epsilon}
\nc{\su}{\widehat{\goth{sl}}_2}
\nc{\sw}{\goth{sl}}
\nc{\g}{\goth{g}}
\nc{\h}{\goth{h}}
\nc{\n}{\goth{n}}
\nc{\N}{\widehat{\n}}
\nc{\ab}{\goth{a}}
\nc{\G}{\widehat{\g}}
\nc{\De}{\Delta_+}
\nc{\gt}{\widetilde{\g}}
\nc{\Ga}{\Gamma}
\nc{\one}{{\bold 1}}
\nc{\hh}{\widehat{\h}}
\nc{\z}{{\goth Z}}
\nc{\zz}{{\cal Z}}
\nc{\Hh}{{\cal H}}
\nc{\qp}{q^{\frac{k}{2}}}
\nc{\qm}{q^{-\frac{k}{2}}}
\nc{\La}{\Lambda}
\nc{\wt}{\widetilde}
\nc{\qn}{\frac{[m]_q^2}{[2m]_q}}
\nc{\cri}{_{\on{cr}}}
\nc{\sun}{\widehat{\sw}_N}
\nc{\HH}{{\cal H}_q(\sw_N)}
\nc{\ca}{\wt{{\cal A}}_{h,k}(\sw_2)}
\nc{\si}{\sigma}
\nc{\gl}{\goth{g}\goth{l}_N}
\nc{\el}{\ell}
\nc{\s}{t}
\nc{\tN}{\theta_{p^N}}
\nc{\ds}{\displaystyle}
\nc{\Dp}{D_{p^{-1}}}
\nc{\Dq}{D_{q^{-1}}}
\nc{\tL}{{\bold L}}
\nc{\tP}{{\bold P}}
\nc{\tA}{{\bold A}}
\nc{\beq}{\begin{equation}}
\nc{\PN}{{\cal P}_{q,N}}
\nc{\RN}{{\cal R}_{q,N}}
\nc{\LN}{{\cal L}_{q,N}}
\nc{\wLN}{\widetilde{\cal L}_{q,N}}
\nc{\MN}{{\cal M}_{q,N}}
\nc{\FN}{{\cal F}_{q,N}}
\nc{\GN}{{\cal G}_{q,N}}
\nc{\JN}{{\cal J}_{q,N}}
\nc{\IN}{{\cal I}_{q,N}}
\nc{\wGN}{\widetilde{\cal G}_{q,N}}
\nc{\Res}{\on{Res}}
\nc{\bi}{\bibitem}

\title[Deformations of soliton equations]{Deformations of the KdV
hierarchy and related soliton equations}

\author{Edward Frenkel}\thanks{Partially supported by NSF grant DMS-9501414}
\address{Department of Mathematics, Harvard University, Cambridge, MA
02138, USA}

\date{November 1995}

\maketitle

\begin{abstract}
We define hierarchies of differential--$q$-difference equations, which
are $q$--deformations of the equations of the generalized KdV
hierarchies. We show that these hierarchies are bihamiltonian, one of
the hamiltonian structures being that of the $q$--deformed classical
${\cal W}$--algebra of $\sw_N$, defined by Reshetikhin and the author.
We also find $q$--deformations of the mKdV hierarchies and the affine
Toda equations.
\end{abstract}

\section{Introduction.}
The $N$th Korteweg-de Vries (KdV) hierarchy is a bihamiltonian integrable
system defined on the space ${\cal M}'_N$ of $N$th order differential
operators on the circle of the form \beq \label{l} \pa^N + u_1(z) \pa^{N-2}
+ \ldots + u_{N-2}(z) \pa + u_{N-1}(z),
\end{equation}
where $\pa \equiv \pa/\pa z$ \cite{A,GD}. The space of local functionals on
${\cal M}'_N$, considered as a Poisson algebra with respect to one of the
hamiltonian structures of this system, is the classical ${\cal
W}$--algebra associated to $\sw_N$, or the Adler--Gelfand--Dickey algebra.

Recently, N.~Reshetikhin and the author \cite{FR} defined a Poisson algebra
$\W_q(\sw_N)$ which is a $q$--deformation of the classical $\W$--algebra of
$\sw_N$. This Poisson algebra consists of functionals on the space ${\cal
M}'_{q,N}$ of $N$th order $q$--difference operators of the form
$$D^N - t_1(z) D^{N-1} + \ldots + (-1)^{N-1} t_{N-1}(z) D + (-1)^N,$$ where
$[D \cdot f](z) = f(zq)$. It is natural to ask whether there is an
integrable hierarchy on ${\cal M}'_{q,N}$, which is hamiltonian with
respect to this Poisson structure. In this work we construct such a
hierarchy and show that it also has another compatible hamiltonian
structure, i.e. is bihamiltonian. Moreover, the hamiltonians of this
$q$--deformed KdV hierarchy have the following heredity property: the
hamiltonian of the $n$th equation of the hierarchy with respect to one of
the Poisson structures coincides with the hamiltonian of the $(n+N)$th
equation with respect to the other structure. This property is well-known
for the hamiltonians of the ordinary KdV hierarchies \cite{A,GD}.

The KdV hierarchies are closely related to the modified KdV (mKdV)
hierarchies and the affine Toda equations, see
\cite{DS,KW,KW1,FF:laws}. Very briefly, the $N$th mKdV hierarchy is the
pull-back of the $N$th KdV hierarchy by the so-called Miura
transformation. The hamiltonian structure of the mKdV hierarchy is that of
a classical Heisenberg algebra. The $N$th affine Toda equation is a
non-local equation on the phase space of the $N$th mKdV hierarchy, whose
local conservation laws are the hamiltonians of the mKdV hierarchy.

Using a $q$--deformation of the Miura transformation defined in \cite{FR},
we construct a $q$--deformations of the $N$th mKdV hierarchy and the affine
Toda equation. We show that our $q$--deformed mKdV hierarchy is
hamiltonian, and its hamiltonians are conserved with respect to the
evolution of the $q$--deformed affine Toda equation.

This paper was greatly inspired by my joint work with B.~Feigin and
A.~Odesskii \cite{FFO} on commutative subalgebras of the $q$--deformed
classical Virasoro algebra (i.e. $\W_q(\sw_2)$) and its elliptic
deformations. The heredity property in the case of $\W_q(\sw_2)$ is a
result of \cite{FFO}, which is generalized here to the case of
$\W_q(\sw_N)$.

In this paper we use the general technique of discrete Lax equations,
developed by B.~Kuperschmidt \cite{Ku}. Although he considers in \cite{Ku}
another class of Lax operators -- difference rather than $q$--difference,
many of his ideas and results can be applied in our context. But there are
also many differences with \cite{Ku}, most importantly, in hamiltonian
formalism. We would also like to mention a recent paper of D.~Gieseker
\cite{Gie}, where another differential-difference deformation of the KdV
hierarchy is considered.

It would be interesting to find solutions of the $q$--deformed hierarchies,
in particular, analogues of the soliton solutions of the ordinary
hierarchies. It would also be interesting to associate such $q$--deformed
hierarchies to other Lie algebras.

The paper is arranged as follows. In Sect.~2 we recall basic facts about
the KdV, mKdV and the affine Toda equations in the form suitable for
generalization. In Sect.~3 we construct the $q$--deformed KdV hierarchies,
and in Sects.~4 and 5 we construct the $q$--deformed mKdV hierarchies and
affine Toda equations, respectively. The proofs of the results presented
here will appear in \cite{F:new}.

\section{The classical equations.}
\subsection{The KdV hierarchy.}
Let us we recall the construction of the $N$th KdV hierarchy. For more
details, the reader may consult the original papers \cite{A,GD} and the
book \cite{D}. Consider the $N$th order differential operator $L$ of the
form \beq \label{L} \pa^N + u_0(z) \pa^{n-1} + u_1(z) \pa^{N-2} + \ldots +
u_{N-2}(z) \pa + u_{N-1}(z).
\end{equation}
Let ${\cal R}_N = \C[\pa^n u_i]_{i=0,\ldots,N-1;n\geq 0}$ be the ring of
differential polynomials in $u_i(z)$'s, and ${\cal P}_N$ be the ring of
pseudo-differential operators with coefficients in ${\cal R}_N$. For any
rational number $\al$, there is a unique $\al$th power $L^\al$ of $L$ in
${\cal P}_N$. The $n$th equation of the $N$th KdV hierarchy can be written
in the Lax form as \beq
\label{lax} \pa_{\tau_n} L = [ L,( L^{n/N} )_+ ],
\end{equation}
where $P_+$ is the differential part of the pseudo-differential operator
$P$. The non-trivial equations correspond to $n$, which are not divisible
by $N$.

\subsection{Hamiltonian structures of the KdV hierarchy.} The equations
\eqref{lax} can be written in hamiltonian form with respect to two
different Poisson structures \cite{A,GD,D}. These two Poisson structures
are compatible in the sense that any linear combination of them is again a
Poisson structure. Moreover, the hamiltonian of the $n$th equation of the
hierarchy with respect to the first structure coincides with the
hamiltonian of the $(n+N)$th equation with respect to the second
structure. In other words, for each $n$ not divisible by $N$ there exists a
local functional of $u_i(z)$'s, $H_n$, such that the $n$th equations can be
written as $$\pa_{\tau_n} L = \{ L,H_{n+N} \}_1$$ and as
$$\pa_{\tau_n} L = \{ L,H_n \}_2,$$ where the subscript distinguishes the
first and the second Poisson structures. This property was discovered by
Magri \cite{Ma} for $N=2$, and by Gelfand and Dickey \cite{GD} for general
$N$. We call it the heredity property.

There are essentially two formulas for the hamiltonians of the KdV
hierarchies. One of them is \cite{A,GD}
$$H_n = \frac{n}{N} \int \on{Res} L^{n/N} dz,$$ where $\on{Res} P$ stands
for the $\pa^{-1}$--coefficient of $P$. The other formula is $H_n = N \int
f_n dz$, where $f_n$ is uniquely determined by the formula
$$\pa = L^{1/N} + \sum_{i\geq 0} f_i L^{-i/N},$$ see \cite{CWF}. From these
formulas we see that $H_n$ is a local functional, i.e. it has the form
$\int R dz$, where $R \in {\cal R}_N$.

V.~Drinfeld and V.~Sokolov \cite{DS} have shown that the phase space ${\cal
M}_N$ of the $N$th KdV hierarchy, which consists of the differential
operators of the form \eqref{L}, can be obtained by hamiltonian reduction
from a hyperplane in the dual space to the affine algebra
$\widehat{\gl}$. The Poisson structure induced on ${\cal M}_N$ by this
reduction coincides with the second Poisson structure of the corresponding
KdV hierarchy. The space of all local functionals on ${\cal M}_N$ is the
classical $\W$--algebra $\W_N$ associated to $\gl$. The hamiltonians $H_n$
of the KdV hierarchy commute with each other with respect to both Poisson
structures. Thus, they span an infinite-dimensional commutative subalgebra
of the $\W$--algebra.

\subsection{Reduction to the submanifold $u_0(z)=0$.} Equations \eqref{lax}
imply: $\pa_{\tau_n} u_0(z) = 0$. Therefore we can set $u_0(z)$ to be
equal to any function. The standard choice is $u_0(z)=0$. Then formula
\eqref{lax} defines a hierarchy of equations on the submanifold ${\cal
M}'_N \subset {\cal M}_N$ of operators of the form \eqref{l}. The two
Poisson structures on ${\cal M}_N$ discussed above can be restricted to
${\cal M}'_N$, and the equations are hamiltonian with respect to these
restrictions. The Poisson algebra of local functionals on ${\cal M}'_N$ is
the classical $\W$--algebra associated to $\sw_N$.

\subsection{The mKdV hierarchy.} Now consider the space ${\cal F}_N$ of
$N$--tuples of first order differential operators
\beq    \label{first}
(\pa+v_1(z),\ldots,\pa+v_N(z)).
\end{equation}
Consider the map $\mu_i: {\cal F}_N \arr {\cal M}_N$, which sends an
$N$--tuple \eqref{first} to $$L_i = (\pa + v_i(z))(\pa + v_{i+1}(z))
\ldots (\pa + v_{i+N-1}(z)).$$ Here and below we identify all indices
modulo $N$. This map is called the $i$th Miura transformation. We can
pull back the KdV hierarchy from ${\cal M}_N$ to ${\cal F}_N$ using
one of these maps. The corresponding hierarchy on ${\cal F}_N$ is
called the $N$th mKdV hierarchy.

The equations of this hierarchy can be written in the Lax form as follows.
Consider the operator \beq
\label{tL} \tL = \begin{pmatrix} 0 & \pa + v_1 & 0 & \hdots & 0 \\ 0 & 0 &
\pa + v_2 & \hdots & 0 \\ \hdotsfor{5} \\ 0 & 0 & 0 & \hdots & \pa +
v_{N-1} \\ \pa + v_N & 0 & 0 & \hdots & 0
\end{pmatrix}.
\end{equation}
There exists a unique matrix pseudo-differential operator of the form $$\tP
= \begin{pmatrix}
P_1 & 0 & \hdots & 0 \\
0 & P_2 & \hdots & 0 \\
\hdotsfor{4} \\
0 & 0 & \hdots & P_N
\end{pmatrix},$$ where $$P_i = \pa + \sum_{m\geq 0} P_i^{[m]} \pa^{-m},$$
and $P_i^{[m]}$ are differential polynomials in $v_i, i=1,\ldots,N$, such
that $[\tL,\tP] = 0$. The $n$th equation of the $N$th mKdV hierarchy reads
\beq \label{tlax} \pa_{t_n} \tL = [ \tL,(\tP^n)_+ ].
\end{equation}
In writing the mKdV equations in this form we followed the idea of
Kuperschmidt \cite{Ku}, Ch.~IV. We could not find this particular
representation of mKdV hierarchy in the literature; another representation
has been given by Kuperschmidt and Wilson \cite{KW1}.

\subsection{Hamiltonian structure of the mKdV hierarchy.}
The space $\F_N$ is isomorphic to the product of a hyperplane in the
dual space to the homogeneous Heisenberg subalgebra of $\sun$ and the
space of functions on the circle. Therefore it has a Poisson
structure, which is the product of the Kirillov-Kostant Poisson
structure on the hyperplane and the trivial structure on the space of
functions. We call the Poisson algebra of local functionals on ${\cal
F}_N$, the classical Heisenberg algebra and denote it by ${\cal H}_N$.
Explicitly, we have:
\begin{align*}
\{ v_i(z),v_i(w) \} &= - \frac{N-1}{N} \delta'\left( \frac{w}{z} \right), \\
\{ v_i(z),v_j(w) \} &= \frac{1}{N} \delta'\left( \frac{w}{z} \right), \quad
i\neq j,
\end{align*}
where $\delta'(x) = \sum_{m\in\Z} m x^m$. In particular, $v_1(z)+
\ldots + v_N(z)$ lies in the kernel of this Poisson structure.

Equation \eqref{tlax} is hamiltonian with respect to this Poisson
structure
\cite{DS,KW,KW1}, i.e. it can be written as $$\pa_{t_n} v_i(z) = \{
v_i(z),{\bold H}_n \}, \quad \quad i=1,\ldots,N,$$ where $${\bold H}_n =
\frac{1}{n} \int \on{Res} \on{tr} \tP^n dz.$$ Moreover, the map $\mu_1$ is
hamiltonian, if we consider ${\cal M}_N$ as a Poisson manifold with respect
to the second Poisson structure \cite{KW1}. Thus, $\mu_1$ defines a
homomorphism of Poisson algebras ${\cal W}_N \arr {\cal H}_N$. The
hamiltonian ${\bold H}_n$ is the image of the $n$th hamiltonian $H_n$ of
the corresponding KdV hierarchy under this homomorphism.

\begin{rem}
Explicitly, $P_i$ is the pull-back of $(L_i)^{1/N}$ under the $i$th Miura
transformation.

Equations \eqref{tlax} imply that $\pa_{t_n} (v_1(z) + \ldots +
v_N(z)) = 0$. Therefore we can set $\sum_{i=1}^N v_i(z) = 0$. This
gives us a submanifold ${\cal F}'_N \subset {\cal M}_N$. The
restriction of $\mu_1$ to ${\cal F}'_N$ is a hamiltonian map ${\cal
F}'_N \arr {\cal M}'_N$.\qed
\end{rem} \medskip

\subsection{The affine Toda equations.} Closely related to the mKdV
hierarchy is the affine Toda equation \beq
\label{toda} \pa_t v_i(z) = e^{\phi_{i+1}(z) -
\phi_i(z)} - e^{\phi_i(z) - \phi_{i-1}(z)}, \quad \quad i=1,\ldots,N,
\end{equation}
where $\pa \phi_i(z) = v_i(z)$. The affine Toda equation can be
represented in the Lax form as follows:
\beq    \label{todalax}
\pa_t \tL = [\tL,\tA\tL^{-1}],
\end{equation}
where $\tA$ is a matrix of the form $$\tA = \begin{pmatrix}
0 & A_1 & 0 & \hdots & 0 \\
0 & 0 & A_2 & \hdots & 0 \\
\hdotsfor{5} \\
0 & 0 & 0 & \hdots & A_{N-1} \\
A_N & 0 & 0 & \hdots & 0
\end{pmatrix}.$$
Again, we could not find this representation of the affine Toda
equation in the literature, but similar representations have been
given in \cite{DS,KW}.

\subsection{Hamiltonian structure of the affine Toda equation.}
The affine Toda equation can be written in hamiltonian form \cite{KW}:
$$\pa_t v_i(z) = \{ v_i(z),{\bold H} \}, \quad \quad i=1,\ldots,N,$$
where $${\bold H} = \sum_{i=1}^N \int e^{\phi_{i+1}(z) - \phi_i(z)}
dz.$$ The hamiltonian ${\bold H}$ is not a local functional, and hence
does not lie in the classical Heisenberg algebra ${\cal H}_N$. Because
of that, the Poisson bracket of ${\bold H}$ with any element of ${\cal
H}_N$ can only be defined as a linear operator acting from ${\cal
H}(\sw_N)$ to another vector space, see \cite{KW} and also
\cite{FF:laws,FF:kdv}. However, we can still define local integrals of
motion of the affine Toda equation as elements in the kernel of the
linear operator $\{ \cdot,{\bold H} \}$
\cite{KW,FF:laws,FF:kdv}. These quantities are conserved with respect to the
evolution defined by the equation \eqref{toda}.

Since the affine Toda equation has Lax representation with the same
operator $\tL$ as the equations of the corresponding mKdV hierarchy, the
hamiltonians of the mKdV hierarchy are integrals of motion of the affine
Toda equation \cite{DS,KW}, i.e. $\{ {\bold H}_n,{\bold H} \} = 0$ for all
$n>0$. The converse is also true: the local functionals which commute with
${\bold H}$ are precisely the linear combinations of ${\bold H}_n$
\cite{FF:laws,FF:kdv}. Therefore this property can be taken as the definition
of
the hamiltonians ${\bold H}_n$, as it was done in
\cite{FF:laws,FF:kdv}.

\section{Deformation of the KdV hierarchy.}
\subsection{The phase space and related rings.}
Consider the space ${\cal M}_{q,N}$ of $q$--difference operators of the
form \beq \label{qL} L = D^N - t_1(z) D^{N-1} + \ldots + (-1)^{N-1}
t_{N-1}(z) D + (-1)^N t_N(z),
\end{equation} where $$t_i(z) = \sum_{m\in\Z} t_i[m] z^{-m}$$ is a Laurent
series for each $i=1,\ldots,N$, and $[D \cdot f](z) = f(zq)$. Let us define
the rings ${\cal R}_{q,N}$ and ${\cal P}_{q,N}$ associated to these
$q$--difference operators, which are analogous to the rings ${\cal R}_N$
and ${\cal P}_N$.

Consider the ring ${\cal R}^{(0)}_{q,N} = \C[t_i(zq^j)]_{i=1,\ldots,N;
j\in\Z}$. The operator $D$ acts naturally on ${\cal R}^{(0)}_{q,N}$. The
ring $\RN$ is the smallest ring containing ${\cal R}^{(0)}_{q,N}$ as a
subring, on which the operator $1+D+\ldots +D^{N-1}$ is invertible. To
construct it, we adjoin to ${\cal R}^{(0)}_{q,N}$ the solutions of the
equations $(1+D+ \ldots+D^{N-1}) f(z) = R(z)$ for all $R(z) \in {\cal
R}^{(0)}_{q,N}$. This gives us a ring ${\cal R}^{(1)}_{q,N}$. Then we
adjoin to ${\cal R}^{(1)}_{q,N}$ the solutions of these equations for all
elements $R$ of ${\cal R}^{(1)}_{q,N}$, which are not in ${\cal
R}^{(0)}_{q,N}$, etc. The inductive limit of the rings $R^{(i)}_{q,N}$,
which are obtained this way, is the ring $\RN$.

If we write $t_i(z) = \sum_m t_i[m] z^{-m}$, then the Fourier coefficient
$R[m]$ of an element $R(z) = \sum_m R[m] z^{-m}$ of $\RN$ is a linear
combination of expressions of the form \beq
\label{qloc} \sum_{m_1+\ldots +m_k=M} c(m_1,\ldots,m_k) t_{i_1}[m_1]
\ldots t_{i_k}[m_k],
\end{equation}
where $c(m_1,\ldots,m_k)$ is a rational function in
$q^{m_1},\ldots,q^{m_k}$, which has no poles for integral values of $m_i$'s
if $q$ is generic. We denote the Fourier coefficient $R[0]$ of $R(z)$ by
$\int R(z)$.

Given an operator of the form \eqref{qL}, we can substitute the
coefficients $t_i[m], i=1,\ldots,N-1; m\in\Z$, into an expression like
\eqref{qloc} and get a number. Therefore Fourier coefficients of elements
of $\RN$ define functionals on the space $\MN$. We call all expressions of
the form \eqref{qloc}, $q$--local functionals. The leading term in the
$(1-q)$--expansion of a $q$--local functional is a local functional. Note
that there are $q$--local functionals, which can not be represented as
Fourier coefficients of elements of $\RN$. We denote the space of all
$q$--local functionals by $\wLN$, and its subspace, spanned by the
$0$th Fourier coefficients of elements of $\RN$ by $\LN$. We have a
map $\int: \RN \arr \LN$, which sends $R(z)$ to $R[0]$.

Now let $\PN$ be the ring consisting of ``pseudo-difference''
operators of the form $$\sum_{n\leq M} R_n(z) D^n,$$ where each
$R_n(z) \in \RN$, and it may be non-zero for any negative $n$.

\begin{lem}
There exists a unique element $P$ of $\PN$, such that $P^N$ is equal
to $L$ given by formula \eqref{qL}.
\end{lem}

We will denote $P^n$ by $L^{n/N}$. Also, for each $R = \sum_n R_n D^n$
from $\RN$, we set: $R_+ = \sum_{n\geq 0} R_n D^n$, $R_- = \sum_{n<0}
R_n D^n$, and $\Res R = R_0$.

\subsection{Lax form of the equations.}
For each $n=1,2,\ldots$, the $n$th equation of the $q$--deformed KdV
is by definition the Lax equation \beq \label{qkdv}
\pa_{\tau_n} L = [ L,(L^{n/N})_+ ] = - [ L,(L^{n/N})_- ].
\end{equation}
{}From the first equality it follows that $\pa_{\tau_n} L$ has the form
$\sum_{n=1}^\infty w_n D^n$, while from the second equality it follows
that it has the form $\sum_{n=-\infty}^{N-1} w_n D^n$. Hence $\pa_{\tau_n}$
defined this way is a difference operator of order $N-1$, and formula
\eqref{qkdv} makes sense.

\begin{rem}
The equation \eqref{qkdv} can be viewed as the integrability condition for
the system $$L \Psi = \la \Psi, \quad \quad \pa_{\tau_n} \Psi = (L^{n/N})_+
\Psi.$$ Note that the first formula looks similar to formulas for the spectra
of transfer-matrices appearing in Bethe ansatz, see \cite{FR} and
references therein.\qed
\end{rem} \medskip

We call the hierarchy of Lax equations \eqref{qkdv} the $N$th $q$--deformed
KdV hierarchy. These equations define flows on the space of operators of
the form \eqref{qL} with $t_i(z)$'s belonging to a sufficiently large class
of functions of real or complex variable. In the present work, we do not
discuss analytic aspects of these flows, but focus on the algebraic
structures, which underlie them. Such approach proved very effective in the
case of the ordinary KdV hierarchies.

{}From the algebraic point of view, formula \eqref{qkdv} defines a
derivation $\pa_{\tau_n}$ on the ring $\RN$ and a linear operator
on the space $\LN$. The latter can clearly be extended to a linear
operator $\wLN \arr \wLN$. The following proposition is analogous to
Proposition 1.7 from \cite{Ku}, Ch.~III.

\begin{prop}    \label{commute}
The derivations $\pa_{\tau_n}$ are non-trivial for all $n$ not divisible by
$N$ and they commute with each other. The $q$--local functionals $\int \Res
L^{n/N}$ are conserved with respect to all $\pa_{\tau_m}$,
i.e. $\pa_{\tau_m} \cdot \int \Res L^{n/N} = 0, \forall n,m$. Moreover,
$\int \Res L^{n/N} \neq 0$, if $n$ is not divisible by $N$.
\end{prop}

\begin{rem}
Denote ${\cal D} = 1-D$. The proof of the fact that the functionals $\int
\Res L^{n/N}$ are conserved with respect to all $\pa_{\tau_m}$ relies on
the property of the map $\int: \RN \arr \LN$ that it vanishes on the
image of ${\cal D}$, i.e. $\int ({\cal D} R) = 0, \forall R \in \RN$.
In general, we may want to allow the coefficients $t_i(z)$ of the
operator \eqref{qL} to belong to a class of functions other than
Laurent power series. Accordingly, we may consider another space of
functionals ${\cal L}'_{q,N}$. In that case, conservation laws of the
equations \eqref{qkdv} are $\phi(\Res L^{n/N})$, where $\phi$ is an
arbitrary map $\RN \arr {\cal L}'_{q,N}$, which vanishes on the image
of ${\cal D}$. For example, in some cases we can use the Jackson
integral as $\phi$.\qed
\end{rem}

\begin{rem}
Consider the equations \beq \label{qkp} \pa_{\tau_n} P = [ P,(P^n)_+],
\quad \quad n=1,2,\ldots,
\end{equation}
where $P$ is an operator of the form $D + \sum_{n<0} p_n D^n$. These
equations define a commuting hierarchy of flows on the space of such
operators. If we impose the condition $P^N = L$, where $L$ is of the form
\eqref{qL}, then this hierarchy is equivalent to the one defined by formula
\eqref{qkdv}. If we do not impose this condition, we obtain a hierarchy,
which is $q$--deformation of the KP hierarchy.\qed
\end{rem} \medskip

\subsection{Hamiltonian structures.}
We define two compatible Poisson structures on the space $\MN$, i.e.
we define two Poisson brackets on the space $\wLN$ of functionals on
$\MN$. It is sufficient to define the Poisson brackets between
$t_i[m]$ and $t_j[n]$, or between the power series $t_i(z)$ and $t_j(w)$.

We set
\beq    \label{p1}
\quad \quad \{ t_i(z), t_j(w) \}_1 =
\end{equation}
$$\begin{cases} \displaystyle
\delta\left( \dfrac{wq^{N-j}}{z} \right) t_N(z) t_{i+j-N}(w) - \delta\left(
\dfrac{w}{zq^{N-i}} \right) t_{i+j-N}(z) t_N(w), & i,j \neq N, i+j \geq N\\
0, & \on{otherwise}
\end{cases}$$
and
\begin{align} \notag
\{ t_i(z),t_j(w) \}_2 = & \sum_{m\in\Z} \left( \frac{w}{z}
\right)^m \frac{(1-q^{im})(1-q^{m(N-j)})}{1-q^{mN}} t_i(z) t_j(w) \\
\label{p2} &+ \sum_{r=1}^{\on{min}(i,N-j)} \delta \left( \frac{wq^r}{z}
\right) t_{i-r}(w) t_{j+r}(z) \\ \notag &- \sum_{r=1}^{\on{min}(i,N-j)}
\delta\left( \frac{w}{zq^{j-i+r}} \right) t_{i-r}(z) t_{j+r}(w), \quad
i\leq j.
\end{align}
In these formulas $\delta(x) = \sum_{m\in\Z} x^m$, and we use the
convention that $t_0(z) \equiv 1$.

\begin{rem}
The parameter $q$ in formula \eqref{p2} corresponds to $q^{-2}$ in notation
of \cite{FR}.\qed
\end{rem} \medskip

These Poisson structures are $q$--deformations of the first and second
Poisson structures of the $N$th KdV hierarchy, respectively. The Poisson
structure $\{\cdot,\cdot\}_2$ given by \eqref{p2} was defined in \cite{FR}
(modulo the relation $t_N(z)=1$). The space $\wLN$, considered as a Poisson
algebra with respect to this structure, is a $q$--deformation of the
classical $\W$--algebra of $\gl$.

Note that these two brackets preserve the subspace $\LN$ of $\wLN$ and
hence define two Poisson brackets on it as well.

\begin{prop}
The Poisson structures $\{ \cdot,\cdot \}_1$ and $\{ \cdot,\cdot \}_2$ are
compatible in the sense that any linear combination of them is also a
Poisson structure.
\end{prop}

The next theorem shows that the $N$th $q$--deformed KdV hierarchy defined
by \eqref{qkdv} is bihamiltonian with respect to these two Poisson
structures. Moreover, the hamiltonians of the hierarchy have the heredity
property characteristic for the hamiltonians of the KdV hierarchies.

\begin{thm}    \label{one}
{\em Let $H_n = \frac{N}{n} \int \Res L^{n/N}$. The equation $\pa_{\tau_n}
L = [L,(L^{n/N})_+]$ can be represented in hamiltonian form as
$$\pa_{\tau_n} L = \{ L,H_{n+N} \}_1 \quad \quad and \quad \quad
\pa_{\tau_n} L = \{ L,H_n \}_2.$$ The hamiltonians $H_n$ commute with each
other with respect to both Poisson structures: $\{ H_n,H_m \}_1 = \{
H_n,H_m \}_2 = 0.$}
\end{thm}

The proof of this theorem will appear in \cite{FFO} for the case $N=2$ and
in \cite{F:new} for the general case.

\begin{rem} The first $N-1$ hamiltonians $H_n, n=1,\ldots,N-1$, commute
with all elements of $\wLN$ with respect to the Poisson structure $\{
\cdot,\cdot \}_1$.\qed
\end{rem} \medskip

\subsection{Reduction to the submanifold $t_N(z)=1$.} Equations
\eqref{qkdv} imply that $\pa_{\tau_n} t_{N(z)} = 0$. Therefore we can
set $t_N(z)=1$. This gives us a submanifold ${\cal M}'_{q,N} \subset {\cal
M}_{N,q}$, which consists of the $q$--difference operators of the form \beq
\label{glN} L = D^N - t_1(z) D^{N-1} + \ldots + (-1)^{N-1} t_{N-1}(z) D +
(-1)^N.
\end{equation}
Formula \eqref{qkdv} defines a hierarchy of differential--$q$-difference
equations on the space ${\cal M}'_{q,N}$. Formulas \eqref{p1} and
\eqref{p2}, in which we set $t_N(z)=1$, define two compatible Poisson
structures on ${\cal M}'_{q,N}$ with respect to which the equations
\eqref{qkdv} are hamiltonian. The Poisson algebra corresponding to the
second structure is a $q$--deformation of the classical $\W$--algebra of
$\sw_N$, defined in \cite{FR}.

There is a nice formula for the first Poisson structure \eqref{p1} on
${\cal M}'_{q,N}$, which is analogous to the formula defining the first
Poisson structure of the ordinary KdV hierarchies. As we mentioned earlier,
elements of the space $\wLN$ can be considered as functionals on the space
$\MN$ of $q$--difference operators of the form \eqref{qL}. On the other
hand, each element $X$ from $\PN$ of the form $\sum_{i=1}^{N-1} x_i D^{-i}$
defines a linear functional $\ell_X$ on $\MN$ by the formula $$\ell_X(L) =
\int \Res LX.$$ The Poisson structure \eqref{p1} is uniquely determined by
the following formula for the Poisson bracket of the linear functionals:
$$\{ \ell_X,\ell_Y \}_1(L) = \int \Res (L[X,Y]).$$ It would be interesting
to find a similar formula for the second Poisson structure.

\subsection{The second construction of hamiltonians.}
Following Kuperschmidt \cite{Ku}, Ch.~IX, we can give another
construction of the hamiltonians $H_n$. The idea of this construction goes
back to \cite{CWF} (see Sect.~2 above).

Let again $P$ be the $N$th root of $L$ given by \eqref{qL}. We can uniquely
represent $D$ as $$D = P + \sum_{i\geq 0} f_i P^{-i},$$ where $f_i \in
\RN$. Now introduce elements $h_n \in \RN$ via the generating series
$$\sum_{n>0} h_n t^n = - \log \left( 1 + \sum_{i\geq 0} f_n t^{-n-1}
\right),$$ where $$- \log (1-x) = \sum_{m>0} \frac{x^m}{m}.$$

\begin{prop}
$h_n$ equals $\frac{1}{n} \Res L^{n/N}$ up to a total difference,
i.e. $$f_n - \frac{1}{n}\Res L^{n/N} = (1-D) g_n$$ for some $g_n \in
\RN$. Therefore $H_n = N \int f_n$.
\end{prop}

The proof of this proposition is analogous to the proof of Theorem 2.24 of
\cite{Ku}, Ch.~IX. Following \cite{Ku} we can define the $\tau$--function
of the $N$th $q$--deformed KdV hierarchy in the same way as for the
ordinary hierarchies.

\subsection{Examples and $q \arr 1$ limit.}
We have: $$P = L^{1/N} = D + b(z) + \ldots,$$ where $b(z)$ is defined by
the equation $$(1+D+\ldots+D^{N-1}) b(z) = - t_1(z).$$

The first equation of the $N$th $q$--deformed KdV hierarchy is therefore:
$$\pa_{\tau_1} L = [L,D+b(z)],$$ which gives \beq
\label{qkdv1} \pa_{\tau_1} t_i(z) = t_i(z)(b(zq^{N-i})-b(z)) + t_{i+1}(zq)
- t_{i+1}(z),
\end{equation}
where $t_{N+1}(z)=0$. The hamiltonian with respect to the second
structure is $$H_1 = N \int \Res P = N b[0] = -t_1[0].$$ From formula
\eqref{p2} we find that equation \eqref{qkdv1} is equivalent to the
equation $$\pa_{\tau_1} t_i(z) = \{ t_i(z),H_1 \}_2.$$

Let us now consider the case $N=2$ more closely. We can reduce the
$q$--deformed KdV equations to the manifold ${\cal M}'_{2,q}$ of
operators of the form $L = D^2 - t(z) D + 1$. To simplify notation let
us write $t_n$ for $t[n]$. The first two non-trivial hamiltonians are:
$H_1 = -t_0$, and $$H_3 = -\frac{1}{2} t_0 + \frac{1}{3}
\sum_{i+j+k=0}
\frac{1}{(1+q^i)(1+q^j)(1+q^k)} t_i t_j t_k.$$ The Poisson brackets are:
$$\{ t(z),t(w) \}_1 = \delta \left( \frac{wq}{z} \right) - \delta
\left( \frac{w}{zq} \right),$$
\beq \label{qvir} \{ t(z),t(w) \}_2 = \sum_{m\in\Z} \left(
\frac{w}{z} \right)^m \frac{1-q^m}{1+q^m} t(z)t(w) + \delta \left(
\frac{wq}{z} \right) - \delta \left( \frac{w}{zq} \right)
\end{equation}
(see \cite{FR}).

The equations of the hierarchy are $\pa_{\tau_n} t(z) = \{ t(z),H_n
\}_2$. For $n=1$ have from \eqref{qkdv1}: $$\pa_{\tau_1} t(z) =
t(z)(b(zq)-b(z)),$$ where $b(z)$ satisfies $b(z)+b(zq)=-t(z)$.
In terms of Fourier coefficients, we have:
\begin{align*}
\pa_{\tau_1} t(z) &= - \sum_{i,j}
\frac{1-q^{i+j}}{(1+q^i)(1+q^j)} t_i t_j z^{-i-j}, \\
\pa_{\tau_3} t(z) &= -
\frac{3}{2} \sum_{i,j} \frac{1-q^{i+j}}{(1+q^i)(1+q^j)} t_i t_j z^{-i-j} \\
&+ \sum_{i,j,k,l}
\frac{1-q^{i+j+k}}{(1+q^i)(1+q^j)(1+q^{-i-j})(1+q^{i+j+k})}
t_i t_j t_k t_l z^{-i-j-k-l}.
\end{align*}

The equations of the $N$th KdV hierarchy can be recovered from the
equations of the $N$th $q$--deformed KdV hierarchy in the limit $q
\arr 1$.  Let us explain this in the case $N=2$. Set $t(z) = 2 - h^2
u(z) + o(h^2)$, where $q=e^h$. Then $$D^2 - t(z)D + 1 \mapsto h^2(\pa^2
+ u(z)) + o(h^2).$$ Hence the manifold ${\cal M}'_{2,q}$ becomes
${\cal M}'_2$. In \cite{FR} it was shown that if we multiply the right
hand side by $h$, then the leading term of the Poisson structure
\eqref{qvir} gives us in the limit $q \arr 1$ the Poisson structure of
the classical Virasoro algebra, which is the second Poisson structure
of the KdV hierarchy. It is also easy to see that the first Poisson
structure of the $q$--deformed KdV hierarchy becomes in the limit $q
\arr 1$ the first Poisson structure of the KdV hierarchy.

The $h$--expansion of the $n$th hamiltonian of the $q$--deformed KdV
hierarchy is $H_n = \on{const} + h^{n+1} H_n^{(0)} + o(h^{n+1})$,
where $H_n^{(0)}$ is the $n$th hamiltonians of the KdV hierarchy.
Therefore if we rescale the derivation $\pa_{\tau_n}$ as
$\pa'_{\tau_n} = h^{n+1} \pa_{\tau_n}$, then the equation
$\pa'_{\tau_n} t(z) = \{ t(z),H_n \}_2$ gives us in the limit $q
\arr 1$ the $n$th equation of the KdV hierarchy. In particular, the
equation corresponding to $H_2$ (see above) gives us the KdV equation
itself.

Note that the hamiltonian $H_n$ is non-zero if and only if $n$ is odd.
The hamiltonians $H_1,H_3,\ldots$, of the $q$--deformed KdV hierarchy
have odd degrees as series in $t_m$. It will be shown in \cite{FFO}
that there exists another infinite set $J_2,J_4,\ldots$, of series in
$t_m$ of even degrees, which commute with each other and with all
$H_n$ with respect to both Poisson brackets. The first of them is
$$J_2 = \sum_{i>0} \frac{i q^i}{1-q^{2i}} t_i t_{-i}.$$ Technically,
the series $J_n$ do not lie in the Poisson algebra ${\cal L}_{2,q}$.
But we can enlarge ${\cal L}_{2,q}$ by adjoining series of the form
\eqref{qloc}, where $c(m_1,\ldots,m_k)$ polynomially depend on
$m_1,\ldots,m_k$. The series $H_n$'s and $J_n$'s will then constitute
a commutative subalgebra in this Poisson algebra. We expect that the
elements $J_n$ also satisfy the heredity property, which the $H_n$'s
satisfy. It would be interesting to construct explicitly the
corresponding hamiltonian equations.

The hamiltonians of the $N$th $q$--deformed KdV hierarchy have degrees
which are not divisible by $N$. We conjecture that there exists an
additional set of series of degrees divisible by $N$, which commute with
these hamiltonians with respect to both Poisson structures.

\begin{rem}
In \cite{Ku} the equations of the Toda lattice hierarchy (which should not
be confused with the affine Toda equations) were represented in the Lax form
as $\pa_n L = [L,(L^n)_+]$, where $L = \zeta + p + v \zeta^{-1}$ and
$\zeta$ is a difference operator. Although this $L$--operator resembles the
$L$--operator of the $q$--deformed KdV hierarchy $D^2 - t_1(z) D + t_2(z)$,
the equations of the two hierarchies are different.\qed
\end{rem}

\section{Deformations of the mKdV hierarchies.}
\subsection{The phase space.}
Consider the space $\FN$ of $N$--tuples of $q$--difference operators
\beq    \label{Ntuple}
(D-\La_1(z),\ldots,D-\La_N(z)).
\end{equation}
We define the ring $\JN$ in the same way as the ring $\RN$. Let
$\JN^{(0)}$ be the ring generated by $\La_i(zq^j), i=1,\ldots,N;
j\in\Z$. We take as $\JN$, the smallest ring containing $\JN^{(0)}$ as
a subring, on which the operator $(1+D+\ldots+D^{N-1})$ is invertible.
Let $\GN$ be the span the $0$th Fourier coefficients of elements of
$\JN$. The space $\GN$ lies in a larger space $\wGN$, which is defined
in the same way as $\wLN$. We have a map $\int: \JN \arr \GN$.

Introduce the ring $\IN$ of ``pseudo-difference'' operators of the
form $\sum_{n\leq M} J_n D^n$, where each $J_n \in \JN$.

In \cite{FR} the $q$--deformed Miura transformation was defined. The $i$th
$q$--deformed Miura transformation $\mu_{i,q}$ is the map $\FN \arr \MN$,
which sends the $N$--tuple \eqref{Ntuple} to \beq    \label{Li} L_i =
(D-\La_i(z))(D-\La_{i+1}(z)) \ldots (D-\La_{N+i-1}(z)),
\end{equation}
where we identify $\La_{N+j}(z) \equiv \La_j(z)$.

\subsection{Lax form of the equations.}
We want to write down the equations on $\FN$, which are the pull-backs of
the equations of the $N$th $q$--deformed KdV hierarchy. The idea of the
following construction of these equations is due to Kuperschmidt, see
\cite{Ku}, Ch.~IV.

Consider the matrix
\beq    \label{qtL}
\tL = \begin{pmatrix}
0 & D-\La_1 & 0 & \hdots & 0 \\
0 & 0 & D-\La_2 & \hdots & 0 \\
\hdotsfor{5} \\
0 & 0 & 0 & \hdots & D-\La_{N-1} \\
D-\La_N & 0 & 0 & \hdots & 0
\end{pmatrix}.
\end{equation}

\begin{lem}
There exists a unique matrix $$\tP =
\begin{pmatrix}
P_1 & 0 & \hdots & 0 \\
0 & P_2 & \hdots & 0 \\
\hdotsfor{4} \\
0 & 0 & \hdots & P_N
\end{pmatrix},$$ where each $P_i$ is an element of $\IN$ of the form
$$P_i = D + \sum_{m\geq 0} P_i^{[m]} D^{-m},$$
which commutes with $\tL$.
\end{lem}

For each $n=1,2,\ldots$, consider the Lax equation \beq \label{qmkdv}
\pa_{t_n} \tL = [ \tL,(\tP^n)_+ ].
\end{equation}
We call the hierarchy of equations \eqref{qmkdv} the $N$th $q$--deformed
mKdV hierarchy.

Since $$\tL^N = \begin{pmatrix} L_1 & 0 & \hdots & 0 \\ 0 & L_2 & \hdots &
0 \\ \hdotsfor{4} \\ 0 & 0 & \hdots & L_N
\end{pmatrix},$$ where $L_i$ is given by formula \eqref{Li}, we see that
$P_i$ is the pull-back of $P = (L_i)^{1/N}$ via the $i$th
Miura transformation, and so equation \eqref{qmkdv} on
$(\La_1(z),\ldots,\La_N(z))$ implies equation \eqref{qkdv} on each
$L_i, i=1,\ldots,N$. Hence the $q$--deformed mKdV hierarchy is really
a pull-back of the $q$--deformed KdV hierarchy.

We consider $\pa_{t_n}$ as a derivation of the ring $\JN$ and as a
linear operator on the space $\GN$. We have the following analogue
of \propref{commute}.

\begin{prop}    \label{mcommute}
The derivations $\pa_{t_n}$ are non-trivial for all $n$ not divisible by
$N$ and they commute with each other. The functionals $\int \on{Tr} \Res
\tP^n$ are conserved with respect to all $\pa_{t_m}$. Moreover, these
functionals are non-zero, if $n$ is not divisible by $N$.
\end{prop}

\subsection{Hamiltonian form of the equations.}
We define a Poisson structure on $\FN$ by the formulas \cite{FR}
\begin{align}    \label{pbg1}
\{ \La_i(z),\La_i(w) \} & = \sum_{m\in\Z} \left( \frac{w}{z} \right)^m
\frac{(1-q^m)(1-q^{m(N-1)})}{1-q^{mN}} \La_i(z) \La_i(w), \\ \label{pbg2}
\{ \La_i(z),\La_j(w) \} & = - \sum_{m\in\Z} \left(
\frac{wq^{N+i-j-1}}{z} \right)^m \frac{(1-q^m)^2}{1-q^{mN}} \La_i(z)
\La_j(w), \quad i<j.
\end{align}

\begin{thm}
{\em Let ${\bold H}_n = \frac{1}{n} \int \on{Tr} \Res \tP^n$. The equation
$\pa_{t_n} \tL = [\tL,(\tP^n)_+]$ can be represented in hamiltonian form as
$$\pa_{t_n} \La_i(z) = \{ \La_i(z),{\bold H}_n \}, \quad \quad
i=1,\ldots,N.$$ The hamiltonians ${\bold H}_n$ commute with each other: $\{
{\bold H}_n,{\bold H}_m \} = 0.$}
\end{thm}

\begin{prop} If we consider the space $\MN$ as a Poisson manifold with
respect to the Poisson structure $\{ \cdot,\cdot \}_2$, then the first
Miura transformation $\mu_{1,q}$ is a hamiltonian map. This map sends $H_n$
to ${\bold H}_n$.
\end{prop}

Equations \eqref{qmkdv} imply that $$\pa_{t_n} (\La_1(z) \La_2(z) \ldots
\La_N(z)) = 0$$ for all $n$. Therefore the $q$--deformed mKdV hierarchy can
be restricted to the submanifold ${\cal F}'_{q,N}$ of $\FN$ defined by the
equation $$\La_1(z) \La_2(z) \ldots \La_N(z) = 1.$$ The Poisson structure
on $\FN$ can also be restricted to ${\cal F}'_{q,N}$. The Miura
transformation $\mu_{1,q}$ provides a hamiltonian map ${\cal F}'_{q,N} \arr
{\cal M}'_{q,N}$.

Finally, we can recover the ordinary mKdV equations from equations
\eqref{qmkdv} in the limit $q \arr 1$, if we set $\La_i(z) = 1 - h v_i(z) +
o(h)$, where $q=e^h$.

\section{Deformation of the affine Toda equations.}
\subsection{Lax form of the equation.}
Consider the Lax equation \beq \label{qtodalax} \pa_t \tL =
[\tL,\tA\tL^{-1}],
\end{equation}
where $\tA$ is a matrix of the form $$\tA = \begin{pmatrix}
0 & A_1(z) & 0 & \hdots & 0 \\
0 & 0 & A_2(z) & \hdots & 0 \\
\hdotsfor{5} \\
0 & 0 & 0 & \hdots & A_{N-1}(z) \\
A_N(z) & 0 & 0 & \hdots & 0
\end{pmatrix}.$$

{}From this equation we find: \beq \label{find} \pa_t \La_i(z) = A_i(z) -
(D-\La_i(z)) A_{i+1}(z) (D-\La_{i+1}(z))^{-1},
\end{equation}
for all $i=1,\ldots,N$. Since
\begin{multline}
(D-\La_i(z)) A_{i+1}(z) (D-\La_{i+1}(z))^{-1} = \\
A_{i+1}(zq) + (A_{i+1}(zq) \La_{i+1}(z) - A_{i+1}(z)
\La_i(z))(D-\La_{i+1})^{-1},
\end{multline}
we obtain from \eqref{find} that \beq
\label{qdiff} A_i(zq) = \La_{i-1}(z) \La_i(z)^{-1} A_i(z),
\end{equation}
and \beq    \label{qtoda}
\pa_t \La_i(z) = A_i(z) - A_{i+1}(zq).
\end{equation}

The function $A_i(z)$ appearing on the right hand side of this equation is
determined by the $q$--difference equations \eqref{qdiff}. In the limit $q
\arr 1$ this equation gives rise to the differential equation satisfied by
the function $e^{\phi_i(z) - \phi_{i-1}(z)}$ appearing on the right hand
side of the affine Toda equation \eqref{toda}. Indeed, if we set $\La_i(z)
= 1 - h v_i(z) + o(h)$, and $A_i(z) = a_i(z) + O(h)$, we obtain from
\eqref{qdiff}:
$$\pa a_i(z) = (v_i(z) - v_{i-1}(z)) a_i(z),$$ which shows that $a_i(z) =
e^{\phi_i(z) - \phi_{i-1}(z)}$. Therefore equation \eqref{qtoda} becomes in
the limit $q \arr 1$ the $N$th affine Toda equations. We call equation
\eqref{qtoda} the $q$--deformed $N$th affine Toda equation. Note that
$\pa_t(\La_1(z) \ldots \La_N(z)) = 0$. Therefore we can set $\La_1(z)
\ldots \La_N(z) = 1$.

\subsection{Hamiltonian form of the equation.}
We now want to represent the $q$--deformed affine Toda equation in
hamiltonian form. The construction of the corresponding Poisson
structure is analogous to the construction of the Poisson structure of
the ordinary affine Toda equation \cite{KW}, see also
\cite{FF:laws,FF:kdv}.

Introduce new power series $Q_i(z)$ as solutions of the $q$--difference
equations \beq    \label{q} Q_i(zq) = \La_i(z) Q_i(z).
\end{equation}
Let $\JN^+$ be the tensor product $\JN \otimes
\C[Q_1(z),Q_1(z)^{-1},\ldots,Q_N(z),Q_N(z)^{-1}]$ where each $Q_i(z)$
satisfies equation \eqref{q}. Note that according to formula \eqref{qdiff}
we can express $A_i(z)$ in terms of $Q_j(z)$'s:
$$A_i(z) = Q_{i-1}(z) Q_i(z)^{-1}.$$ Let $\GN^+$ be the span of the $0$th
Fourier coefficients of elements of $\JN^+$, and $\int$ be the
corresponding map $\JN^+ \arr \GN^+$.

We can extend the Poisson bracket on $\GN$ to a partial Poisson bracket $\{
\cdot,\cdot \}: \GN \times \GN^+ \arr \GN^+$ by postulating
\begin{align*}
\{ \La_i(z),Q_i(w) \} &= - \sum_{m\in\Z} \left( \frac{w}{z} \right)^m
\frac{1-q^{m(N-1)}}{1-q^{mN}} \La_i(z) Q_i(w), \\
\{ \La_i(z),Q_j(w) \} &= \sum_{m\in\Z} \left(
\frac{wq^{N+i-j-1}}{z} \right)^m \frac{1-q^m}{1-q^{mN}} \La_i(z)
\La_j(w), \quad i<j, \\
\{ \La_i(z),Q_j(w) \} &= \sum_{m\in\Z} \left(
\frac{wq^{i-j-1}}{z} \right)^m \frac{1-q^m}{1-q^{mN}} \La_i(z)
\La_j(w), \quad i>j.
\end{align*}
Note that these formulas are consistent with formulas \eqref{pbg1} and
\eqref{pbg2} if we take into account equation \eqref{q}.

Denote $$S_i(z) = Q_i(w) Q_{i+1}(wq)^{-1}, \quad \quad i=1,\ldots,N.$$
{}From the Poisson brackets above we derive:
\begin{align*}
\{ \La_i(z),S_i(w) \} &= - \delta \left( \frac{w}{z} \right) \La_i(z)
S_i(w) = - \delta \left( \frac{w}{z} \right) A_{i+1}(zq), \\ \{
\La_i(z),S_{i-1}(w) \} &= \delta \left( \frac{w}{z} \right) \La_i(z)
S_{i-1}(w) = \delta \left( \frac{w}{z} \right) A_i(z), \\ \{
\La_i(z),S_j(w) \} &= 0, \quad \quad j\neq i,i-1.
\end{align*}
These formulas lead us to the following result.

\begin{prop}
The $q$--deformed affine Toda equation can be presented in hamiltonian form
as $$\pa_t \La_i(z) = \{ \La_i(z),{\bold H} \}, \quad \quad i=1,\ldots,N,$$
where $${\bold H} = \sum_{i=1}^N \int Q_i(w) Q_{i+1}(wq)^{-1}.$$
\end{prop}

Consider the case $N=2$. We can set $\La_2(z) = \La_1(z)^{-1}$.
Denote $\La(z) = \La_1(z)$. Let $Q(z)$ be a solution of the
equation $Q(zq) = \La(z) Q(z)$. The equation \eqref{qtoda} reads in
this case $$\pa_t \La(z) = Q(z)^{-2} - Q(zq)^2.$$ This is a
$q$--deformation of the sine-Gordon equation. The hamiltonian is
$${\bold H} = Q(z) Q(zq) + Q(z)^{-1}Q(zq)^{-1}.$$

In general, the operator $\pa_t$ defines operators $\JN \arr \JN^+$
and $\GN \arr
\GN^+$. Since the $q$--deformed affine Toda equation can be represented in
the Lax form with the same Lax operator $\tL$ as the equations of the
$q$--deformed mKdV hierarchy, we conclude that the conservations laws
${\bold H}_n$ of this hierarchy are also conserved with respect to the
$q$--deformed affine Toda equation.

\begin{thm} {\em The hamiltonians ${\bold H}_n$ of the $q$--deformed mKdV
hierarchy commute with the hamiltonian of the $q$--deformed affine Toda
equation: $\{ {\bold H}_n,{\bold H} \} = 0$.}
\end{thm}

\begin{rem}
Each summand $\{ \cdot,S_i(z) \}$ of the linear operator $\{ \cdot,{\bold
H} \}: \JN \arr \JN^+$ gives us a map $\JN \arr \JN^i$, where $\JN^i = \JN
\otimes S_i(z)$. Hence $\{ J,{\bold H} \} = 0$ if and only if $\{ J,\int
S_i(z) \} = 0$ for all $i=1,\ldots,N$.\qed
\end{rem}

\subsection{Classical $\W$--algebra as algebra of integrals of motion.}
The hamiltonian equation corresponding to the hamiltonian $\sum_{i=1}^{N-1}
\int e^{\phi_{i+1}(z) - \phi_i(z)} dz$ is the Toda field equation
associated to the Lie algebra $\sw_N$ (finite Toda equation for
shorthand). We can construct a $q$--analogue of this equation:
$$\pa_t \La_i(z) = \{ \La_i(z),\sum_{i=1}^{N-1} \int S_i(z) \}.$$ The
next proposition shows that the algebra of local integrals of motion of
this equation contains the $q$--deformation of the classical $\W$--algebra
of $\sw_N$.

\begin{prop}    \label{kernel}
For all $J \in \JN$, which lie in the image of $\RN$ under the map
$\mu_{1,q}$, $[J,\int S_i(z)] = 0, i=1,\ldots,N-1$.

For all $G \in \GN$, which lie in the image of $\LN$ under the map
$\mu_{1,q}$, $[G,\int S_i(z)] = 0, i=1,\ldots,N-1$.
\end{prop}

We conjecture that the images of $\RN$ and $\LN$ exhaust all elements in
the intersection of kernels of the operators $[\cdot,\int S_i(z)]: \JN \arr
\JN^+$ and $\GN \arr \GN^+$.

The operators $\int S_i(z), i=1,\ldots,N-1$, are the classical limits of
the screening operators constructed in \cite{LP,SKAO} in the case $N=2$ and
in \cite{FF,AKOS} in general case. In \cite{SKAO,FF,AKOS} it was proved
that the quantum $\W$--algebra of $\sw_N$ commutes with the screening
operators. \propref{kernel} is the classical version of this fact.

We can define spaces of quantum integrals of motion of the $q$--deformed
finite and affine Toda equations as the intersections of kernels of the
screening operators corresponding to $\int S_i(z)$, where $i=1,\ldots,N-1$
in the finite case, and $i=1,\ldots,N$ in the affine case. The latter is
clearly a subalgebra of the former. We expect that the integrals of motion
of the quantum affine Toda equation constitute a commutative subalgebra in
the quantum $\W$--algebra. We also expect that the spaces of classical and
quantum integrals of motion are isomorphic as vector spaces. This is true
for the ordinary (i.e. not $q$--deformed) affine Toda equations
\cite{FF:toda,FF:laws}. In that case the quantum integrals of motion can be
interpreted as conservation laws of certain perturbations of conformal
field theories with $\W$--symmtery \cite{int}.

\vspace{5mm}
\noindent{\bf Acknowledgements.} I would like to thank B.~Feigin and
A.~Odesskii for their collaboration in \cite{FFO} which has led to
this work.

\end{document}